\documentclass[tikz,border=5pt]{article}
\usepackage{etoolbox}
\usepackage{amssymb}
\usepackage{amsmath,amsthm}
\usepackage{mathtools}
\usepackage{bm}
\usepackage{extarrows}
\usepackage{indentfirst}
\usepackage[english]{babel}
\usepackage{tikz}
\usepackage{float}
\usepackage{lmodern} 
\usetikzlibrary{calc}
\usetikzlibrary{positioning,shapes,arrows.meta, shapes.geometric,fit}

\usepackage{authblk}

\theoremstyle{thmstyleone}%
\newtheorem{theorem}{Theorem}[section]
\newtheorem{proposition}[theorem]{Proposition}%

\theoremstyle{thmstyletwo}%

\theoremstyle{thmstylethree}%
\newtheorem{definition}{Definition}%
\newtheorem{lemma}{Lemma}

\newtheorem{thm}{Theorem}[section]

\newtheorem{conj}[thm]{Conjecture}

\newtheorem{claim}{Claim}[section]

\newtheorem{constr}{Construction}

\def\ex{\mbox{ex}}

\title{The generalized Tur\'an number for $K_3$ in graphs without suspensions of a path on five vertices
}
\author[1]{Doudou Hei}
\author[2,3]{Xinmin Hou}
\author[4]{Yue Ma}

\affil[1]{School of Mathematical Sciences, Suzhou University of Science and Technology, Suzhou, Jiangsu 215009, China}
\affil[2]{School of Mathematical Sciences, University of Science and Technology of China,  Hefei, Anhui 230026, China}
\affil[3]{Hefei National Laboratory, University of Science and Technology of China,  Hefei, Anhui 230088, China}
\affil[4]{School of Mathematics and Statistics, Nanjing University of Science and Technology,  Nanjing, Jiangsu 210094, China}

\begin{document}

\maketitle

\begin{abstract}
  Given graphs $H$ and $F$, the generalized Tur\'an number $\ex(n, H, F)$ is defined as 
  the maximum number of copies of $H$ in an $n$-vertex graph that contains no copy of $F$.
  The suspension $\widehat{F}$ of a graph $F$ is 
  obtained  by adding a new vertex that is adjacent to every vertex of $F$.
  Mubayi and Mukherjee (2023, DM) conjectured that $\ex(n, K_3, \widehat{P_k})=\left\lfloor \frac{k-2}{2}\right\rfloor \cdot \frac{n^2}{8}+o(n^2)$, where $P_k$ is a
  path on $k\ge 4$ vertices. Using the triangle removal
  lemma, they verified this conjecture for $k=4,5,6$. 
  Later, Mukherjee (2024, DM) established the exact value  $\ex(n, K_3, \widehat{P_4})=\left\lfloor n^2/8\right\rfloor$.
  In this paper, using the stability method, we determine the exact value of $\ex(n, K_3, \widehat{P_5})$ by showing that for sufficiently large $n$, $\ex(n,K_3, \widehat{P_5})=\left\lfloor n^2/8\right\rfloor.$
\end{abstract}

\section{Introduction}\label{sec1}

Let $P_l$, $C_l$
 $K_l$, and $M_l$ denote a path, cycle, complete and almost perfect matching graph on $l$ vertices, respectively. Fix a graph $F$, we say that a graph $G$ is $F$-free if it does not contain $F$ as a subgraph.  Fix graphs $H$ and $G$, we
denote the number of copies (not necessarily induced) of $H$ in $G$ by $N(H, G)$. 
For convenience, let $t(G):=N(K_3, G)$.
The {\em generalized Turán number} is defined as 
$$\ex(n, H , F)=\max\{N(H, G): G \text{ is } F\text{-free}, \vert V(G)\vert=n\}.$$
When $H= K_2$, this is the Turán number $\ex(n, F)$ of graph $F$.
After  decades of isolated results, e.g., \cite{bollobas2008pentagons, erdos1988problems, gyori2012maximum}, the systematic study of $\ex(n, H , F)$ for $H \neq K_2$ was initiated by Alon and Shikhelman~\cite{871}.
Since their work, a tremendous amount of work, e.g., \cite{yuan2025generalized, zhu2023maximum}, has been
done on this function, known as the generalized Tur\'an problems.
See \cite{872} for a comprehensive survey.

The {\emph suspension}  of a graph $F$, denoted by $\widehat{F}$, is the graph obtained by adding an additional vertex to $F$ and connecting it to every vertex of $F$.
Mubayi and Mukherjee \cite{7311} studied $\ex(n, K_3, \widehat{F})$ for various different bipartite graphs $F$. In particular, they investigated $\ex(n, K_3, \widehat{P}_k)$ and gave the following bound and open probelm.
\begin{proposition}[Mubayi and Mukherjee \cite{7311}]\label{lemma 1}
  Let $n\geq k\geq 4$. Then 
  $$\left\lfloor \frac{k-2}{2}\right\rfloor \cdot \frac{n^2}{8}\leq ex(n,K_3, \widehat{P_k})\leq \frac{k-2}{12}\cdot n^2+\frac{(k-2)^2}{12}\cdot n, $$
  where the lower bound holds when $n$ is a multiple of $4\left\lfloor \frac{k-2}{2}\right\rfloor$.
\end{proposition}

The lower bound construction is given by the following graph.

\begin{constr}[Mubayi and Mukherjee \cite{7311}]\label{CONS:Fnk} Let $F_{n, k}$ be the graph formed by  the complete bipartite graph $K_{n/2, n/2}$ with partition $(A, B)$,  together with additional edges in $A$ such that $F_{n, k}[A]$ consists of disjoint copies of $K_{\left\lfloor\frac{k-2}{2}\right\rfloor,\left\lfloor\frac{k-2}{2}\right\rfloor}$,  where $n$ is a multiple of $4\left\lfloor \frac{k-2}{2}\right\rfloor$.
\end{constr}	
From Construction~\ref{CONS:Fnk}, we have 
$$
e(A)=\left\lfloor\frac{k-2}{2}\right\rfloor^2 \cdot \frac{n}{4\left\lfloor\frac{k-2}{2}\right\rfloor}=\left\lfloor\frac{k-2}{2}\right\rfloor \cdot \frac{n}{4}.
$$
Since every triangle of $F_{n, k}$ consists of an edge from $F_{n, k}[A]$ and a vertex from $B$, we have 
$$
t\left(F_{n, k}\right)=\left\lfloor\frac{k-2}{2}\right\rfloor \cdot \frac{n}{4} \cdot \vert A_2\vert=\left\lfloor\frac{k-2}{2}\right\rfloor \cdot \frac{n^2}{8}.
$$
Further, $F_{n, k}$ is $\widehat{P}_k$-free, since the neighborhood of every vertex in $B$ is a disjoint union of $K_{\left\lfloor\frac{k-2}{2}\right\rfloor,\left\lfloor\frac{k-2}{2}\right\rfloor}$, and the neighborhood of every vertex in $A$ is isomorphic to $K_{\left\lfloor\frac{k-2}{2}\right\rfloor, \vert A_2\vert}$.
\\

Mubayi and Mukherjee \cite{7311} believe that the lower bound above is asympototically
tight for all fixed $k\geq 4$ and propose the following conjecture.\\

\begin{conj}[Mubayi and Mukherjee \cite{7311}]\label{Conj:MM23}
Let $n\geq k\geq 4$. Then 
  $$\ex(n,K_3, \widehat{P_k})=\left\lfloor \frac{k-2}{2}\right\rfloor \cdot \frac{n^2}{8}+o(n^2).$$
\end{conj}

They proved the conjecture for the 
first three cases: $k=4, 5$, and 6. 
\begin{theorem}[Mubayi and Mukherjee \cite{7311}]\label{thm 1}
  For $k=4, 5$ and 6, 
  $$\ex(n,K_3, \widehat{P_k})=\left\lfloor \frac{k-2}{2}\right\rfloor \cdot \frac{n^2}{8}+o(n^2).$$
  When $k=4$ or $k=6$, the error term can be improved to $O(n)$.
\end{theorem}

An exact result of $\ex(n, K_3, \widehat{P_4})$ for sufficiently large $n$ was given by Gerbner \cite{873} using the technique of progressive
induction. In particular,  he proved that for a number $K\leq 1575$ and $n\geq 525 + 4K$,
$$\ex(n,K_3, \widehat{P_4})=\left\lfloor n^2/8\right\rfloor.$$ 
Mukheherjee \cite{874}  determined the exact value of
$\ex(n, K_3, \widehat{P_4})$ for every $n\geq 4$, thus closing the gap in the literature for this extremal problem. 
Their method used  induction along with computer programming to prove a base case of the
induction hypothesis.
\begin{theorem}[Mukheherjee~\cite{874}]
  For $n\geq 8$, $\ex(n,K_3, \widehat{P_4})=\left\lfloor n^2/8\right\rfloor$. For $n = 4, 5, 6, 7$, the values of $\ex(n,K_3, \widehat{P_4})$ are
$4, 4, 5, 8$, respectively.
\end{theorem}

The extremal graph considered in \cite{873,874} (for $n \geq 8$) was the following graph.
\begin{constr}
Let $H_n$ be the graph constructed by  adding  a perfect matching to one part of an almost balanced complete bipartite graph on $n$ vertices. 
Specifically:
\begin{itemize}
	\item If $n=4 k$, then $H_n$ is obtained from $K_{2 k, 2 k}$ by adding a perfect matching to either part set.
	
	\item If $n=4 k+1$, then $H_n$ is obtained from  $K_{2 k, 2 k+1}$ by adding a perfect matching to the smaller part set.
	
	\item If $n=4 k+2$, then $H_n$ is obtained from $K_{2 k, 2 k+2}$ by adding a perfect matching to either part set.
	
	\item	If $n=4 k+3$, then $H_n$ is obtained from  $K_{2 k+1,2 k+2}$ by adding a perfect matching to the larger part set.
\end{itemize} 
\end{constr}

Clearly, the neighborhood of every vertex in $H_n$ is either a star or a matching. Thus, $H_n$ is $\widehat{P}_4$-free. 
{A short case analysis shows that the total number of triangles in these graphs is given by $\left\lfloor n^2 / 8\right\rfloor$.}

In this article, using the stability method, we determine the exact value of $\ex(n,K_3, \widehat{P_5})$ for sufficiently large $n$.
Crucially, we verify that this method maintains its validity when applied to $\ex(n,K_3, \widehat{P_4})$. The following is our main result.

\begin{theorem}\label{main 1}
 $\ex(n,K_3, \widehat{P_5})=\left\lfloor n^2/8\right\rfloor,$ where $n$ is sufficiently large.
  
\end{theorem}

It can be verified that $H_n$ is also $\widehat{P_5}$-free, thus
proving the lower bound in Theorem \ref{main 1} for general $n$. We will show that it is in fact  the unique
extremal graph for $\ex(n, K_3, \widehat{P_5})$.



It is worth noting that while  stability methods have been extensively employed in proofs concerning  generalized Tur\'an problems $\ex(n, H, F)$, e.g.,  \cite{gerbner2023generalized, gerbner2023stability}, 
a common feature among these results is that $\chi(F)>\chi(H)$.
By contrast,  the specific problem we address satisfies $\chi(F)=\chi(H)=3$. 
This is the main contribution of our work.

The rest of this article is arranged as follows. In Section 2, we present the preliminary results used in our proofs. The proof of Theorem
\ref{main 1} is then provided in Section 3.

\section{Preliminary}\label{sec2}

In this section, we give preparatory lemmas for the proof of Theorem \ref{main 1}.

First, we need a result on enumerating the edges of extremal configurations which 
can be indirectly derived from the  proof of Theorem 1.4 for $k=4$ in \cite{7311}.

\begin{theorem}[Mubayi and Mukherjee \cite{7311}]\label{main 2}
  Let $G$ be a $\widehat{P_5}$-free  graph on $n$-vertex with $t(G)\geq \frac{n^2}{8}-o(n^2)$,
  then $e(G)\geq \frac{n^2}{4}-o(n^2).$
\end{theorem}


Next, we employ the following theorem to carry out the preliminary characterization of extremal graphs.
Let $T(n, k)$ denote the Tur\'an graph on $n$ vertices with $k$ partites.
\begin{theorem}[Erdős-Simonovits Stability Theorem \cite{821}]\label{thm89}
For $\varepsilon>0$ and any graph $F$ with $\chi(F) \geq 3$, there exists $\delta>0$ such that if $G$ is an $n$-vertex $F$-free graph with

$$
e(G)>\operatorname{ex}(n, F)-\delta n^2,
$$
then the edit distance between $G$ and $T(n, \chi(F)-1)$ is at most $\varepsilon n^2$. In other words, we can add and delete at most $\varepsilon n^2$ edges of $G$ to obtain $T(n, \chi(F)-1)$.
  
\end{theorem}


We need the following classical Tur\'an number of $\widehat{P_5}$~\cite{5}. Let
$$
f(n, k)=\max \left\{n_0 n_1+\left\lfloor\frac{(k-1) n_0}{2}\right\rfloor: n_0+n_1=n\right\}.
$$
Based on the Erdős-Sós Conjecture (\textbf{Erd\H{o}s-S\'os Conjecture:} For any tree $T$, $\ex(n, T)\leq\frac{\vert T\vert-2 }{2}n$.), Zhu, Wang, Zhang, and Zhang~\cite{5} proved the following theorem.
\begin{theorem}[Zhu, Wang, Zhang, and Zhang \cite{5}]
  Let $T$ be a balanced tree of size $2 k$ or $2 k+1$ and Erdős-Sós Conjecture holds for all of its subtrees. When $n \geq 4(4 k)^6$, we have
$$
\ex(n, \widehat{T}) \leq f(n, k)
$$
Moreover, the equality holds for infinitely many $n$.
\end{theorem}

The Tur\'an number for path is given by Erdős and Gallai~\cite{875}.
\begin{theorem}[Erdős and Gallai \cite{875}]\label{thm8121}
  $$\ex(n, P_k)\leq \frac{n}{k-1}\binom{k-1}{2}\leq \frac{k-2}{2}n. $$
\end{theorem}

Obviously, Erdős-Sós Conjecture holds for all paths, by simple calculation, we have

\noindent \textbf{Corollary.}  When $n \geq 4(4 k)^6$, we have
$$\ex(n, \widehat{P_5})\leq\frac{n^2}{4}+\left\lfloor \frac{n+1}{4}\right\rfloor.$$

Next we give some definitions and an important lemma which are instrumental in the proof of Theorem \ref{main 1}.

\begin{definition}
  Let $G$ be a graph with vertex partition $V_1\cup A_2$. Triangles in $G$ are classified into 
  three types: the first type consists of triangles that intersect with $V_1$ at two vertices and with $A_2$ at one vertex; 
  the second type comprises triangles that intersect with $A_2$ at two vertices and with $V_1$ at one vertex; 
  the third type includes triangles whose vertices entirely lie within either $V_1$ or $A_2$.
  For a vertex $u$, let $T_i(u)$ denote the set of triangles of type $i$ in $G$ which contain $u$, where $i=1,2,3$.
  And let $t_i(u)=\vert T_i(u)\vert$.
  For an edge $uv\in E(G)$, let $T_i(uv)$ denote the set of triangles of type $i$ in $G$ which contain $u$ and $v$, where $i=1,2,3$.
  And let $t_i(uv)=\vert T_i(uv)\vert$.

  Let $t(u)$ ($t(uv)$) denote the number of triangles in $G$ contain $u$ ($uv$). Then 
  $t(u)=t_1(u)+t_2(u)+t_3(u)$, $t(uv)=t_1(uv)+t_2(uv)+t_3(uv)$.

  For subgraph $H\subseteq G$, let $t(H)$ denote the number of triangles in $G$ which 
  contain at least one vertex of $H$.

  Let $N_i(u)$ denote the neighborhood of $u$ in $V_i$, and $d_i(u):=\vert N_i(u)\vert$, where $i=1,2$.

\end{definition}

\begin{lemma}\label{8291}
   Let $G$ be a $\widehat{P_5}$-free graph with vertex partition $V_1\cup A_2$, where $\vert V_1\vert=2m$,
   $\vert A_2\vert=k=o(m)$, and $G[V_1]\cong M_{2m}$. If for any vertex $u\in A_2$, 
   the number of non-neighbors of $u$ in $V_1$ is $\Omega(m)$, then 
   $t(G)<k\cdot m.$
\end{lemma}
\begin{proof}
  If there does not exist edge in $G[A_2]$, obviously, $t(G)<k\cdot m.$
  May assume $u,v\in A_2$ and $u\sim v$. Since $G[V_1]\cong M_{2m}$,
  $G[N_1(u)]$ is composed of some isolated vertices and independent edges.
  Let $I(u)$ denote the set of isolated vertices in $N_1(u)$, 
  and $M(u)$ denote the set of vertices that belong to some independent edge.
  Then $N_1(u)=I(u)\cup M(u)$.

  We claim $v$ is adjacent to at most one edge of $G[M(u)]$, namely at most two vertices of $M(u)$, otherwise assume
  $xy,zw\in G[M(u)]$, and $x, z\sim v$, then $yxvzw$ is a copy of 
  $P_5$ in $N(v)$, a contradiction. Then there are at most 2 triangles of 
  $T_2(uv)$ which are not contained in $\{u,v\}\cup I(u)$.

  If there exsits vertex $w\in N_2(u)\setminus\{v\}$, such that 
  $N_1(u)\cap N_1(v)\cap N_1(w)\neq \emptyset$, then $\vert N_1(u)\cap N_1(v)\cap N_1(w) \vert\leq 2$ ( otherwise
  we can find a copy of $\widehat{P_5}$ with center $u$).
  If $\vert N_1(u)\cap N_1(v)\cap N_1(w) \vert= 2$, then neither $v$ nor $w$
  can have other neighbors in $N_1(u)$. And $t_2(uv)=t_2(uw)=2$.
  If $\vert N_1(u)\cap N_1(v)\cap N_1(w) \vert= 1$, then there is at least one vertex
  of $\{v, w\}$ such that it does not have other neighbors in $N_1(u)$. Namely, 
  $t_2(uv)=1$ or $t_2(uw)=1$.

  \begin{claim}
  For any vertex $u\in A_2$, $t_2(u)\leq \vert I(u)\vert +4d_2(u)$.
  \end{claim}
\begin{proof}
  First, we count the number of  triangles in $T_2(u)$ which contain a vertex $v\in N_2(u)$
  with $t_2(uv)\leq 2$. The number is at most $2d_2(u)$.

  Let $X$ be subset of $N_2(u)$ the vertices of which have not been involved before, then 
  $N_1(u)\cap N_1(w)\cap N_1(h)= \emptyset$, for any $w, h\in X$. Since for any vertex $x\in X$, there are at most 2 triangles of 
  $T_2(ux)$ which are not contained in $\{u,x\}\cup I(u)$, we have
  $$\sum_{x\in X}t_2(ux)\leq \vert I(u)\vert +2d_2(u). $$
  Combined with $2d_2(u)$, we have $t_2(u)\leq \vert I(u)\vert +4d_2(u)$.
\end{proof}

Obviously, for any vertex $u\in A_2$, $t_1(u)\leq \frac{\vert M(u)\vert}{2}$. Then we have
\[
\begin{aligned}
t(G)&=\sum_{u\in A_2}t_1(u)+\frac{1}{2}\sum_{u\in A_2}t_2(u)+t(G[A_2]) \\
&\leq \sum_{u\in A_2}\frac{\vert M(u)\vert}{2}+\sum_{u\in A_2} (\frac{\vert I(u)\vert}{2} +2d_2(u))+t(G[A_2])\\
&=\sum_{u\in A_2}\frac{\vert M(u)\vert+\vert I(u)\vert}{2}+\sum_{u\in A_2} 2d_2(u)+t(G[A_2])\\
&\leq km-\Omega(km)+O(k^2)<km,
\end{aligned}
\]
where the penultimate inequality holds since $G[A_2]$ is  $\widehat{P_5}$-free.
\end{proof}

\section{Proof of Theorem \ref{main 1}}

  
Let $G$ be an $n$-vertex $\widehat{P_5}$-free graph with $\mathcal{N}(K_3, G)=\operatorname{ex}(n, K_3, \widehat{P_5})$. 
{Then $\mathcal{N}(K_3, G)\ge \mathcal{N}(K_3, H_n)=\lfloor\frac{n^2}8\rfloor$.} By Theorem \ref{main 2}, $e(G)\geq \frac{n^2}{4}-o(n^2)$.
Since $\ex(n, \widehat{P_5})\leq\frac{n^2}{4}+\left\lfloor \frac{n+1}{4}\right\rfloor$,
by Theorem \ref{thm89}, $G$ can be obtained from a complete bipartite graph $T$ with parts $V_1$ and  $A_2$ 
by adding and removing $o(n^2)$ edges. We choose  $T$ so as to minimize the number of edges that need to be added  and removed 
in this process. 
In particular, every vertex $v \in V_i$ is adjacent to at least as many vertices in $V_{3-i}$ as in $V_i$ 
(otherwise, move $v$ to $V_{3-i}$). Moreover, we have $\vert V_i\vert=\vert A_2\vert-o(n)$, for $i=1,2$.

We may assume 
\begin{equation}\label{EQ:t(x)}
	t(x)\geq \frac{n-6}{4}, \text{ for all $x\in V(G)$.}
\end{equation} 
Indeed, we can assume $n\geq n_0+\binom{n_0}{3}$
for some sufficiently large $n_0$.  If  there exists a vertex
$v_n\in V(G)$ satisfying 
$t(v_n) < \frac{n-6}{4}$, set $G_n:=G$ and $G_{n-1}:=G_n-v_n$, then we have
$$
t(G_{n-1}) =t(G_n)-t(v_n)> t(H_n)-\frac{n-6}{4} \geq t(H_{n-1})+1.
$$
Assume that $G_{\ell}$ on $\ell$ vertices with
$$
t(G_{\ell}) \geq t(H_\ell)+n-\ell
$$
has been defined for some $\ell \leq n-1$. If there exists some vertex $v_{\ell} \in V(G_\ell)$ satisfying $t(v_\ell)< \frac{l-6}{4}$, set $G_{\ell-1}:=G_{\ell}-v_{\ell}$. 
Then we get
$$
t(G_{\ell-1}) =t(G_\ell)-t(v_\ell)> t(H_\ell)+n-\ell-\frac{l-6}{4} \geq t(H_{\ell-1})+n-\ell+1.
$$
Otherwise, terminate the procedure. Let $G_s$ be the final graph when the above iteration terminates. 
Then $G_s$ has exactly $s$ vertices and $t(x)\geq \frac{s-6}{4}$ for all $x\in V(G_s)$. If $s<n_0$, then we have
$$
\binom{n_0}{3}>\binom{s}{3} \geq t(G_s) \geq t(H_s)+n-s \geq n-s>n-n_0 \geq \binom{n_0}{3},
$$
a contradiction. Therefore, we have a subgraph $G_s$ of sufficiently large order $s\left(\geq n_0\right)$ with
 $t(G_s) \geq t(H_s)+n-s$ and $t(x)\geq \frac{s-6}{4}$ for all $x\in V(G_s)$. 
If we can prove $G_s\cong H_s$, then
$$
t(H_s)+n-s \leq t(G_s) \leq t(H_s).
$$
Thus we have $n=s$ and $G_s=G \cong H_n$. 
Therefore, since $s$ is large enough, we can do the same analysis on $G_s$ as $G$. For the sake of writing convenience, in the following proof, we still use $G$ to denote $G_s$.

Since $G$ is $\widehat{P_5}$-free, $G[N(x)]$ is $P_5$-free for all $x\in V(G)$. Thus
$$\frac{3}{2}d(x)\geq e(G[N(x)])=t(x)\geq \frac{n-6}{4},$$ 
where the first inequality holds according to Theorem~\ref{thm8121}.
Therefore, we have $d(x)\geq \frac{n-6}{6}$ for all $x\in V(G)$.

Let $r(u)$ denote the number of edges incident to $u$ in $T$ that are not in $G$, i.e. 
the missing edges between $u$ and vertices in the other part. 
Then we have $\sum_{u \in V(G)} r(u)=o\left(n^2\right)$. 
Thus there are $o(n)$ vertices $u$ with $r(u)=\Omega(n)$. 
Let $A$ denote the set of vertices with $r(u)=o(n)$ and $A_i=A \cap V_i$. 
Then $\left|A_i\right|=\left|V_i\right|-o(n)=\vert A_2\vert-o(n)$.

Let $B_i=V_i\setminus A_i$ for $i=1,2$. Then $\vert B_i\vert=o(n)$. Let $B=B_1\cup B_2$. Then the number of triangles that contain at least one vertex from $B$ is bounded by $o(n^2)$.
In addition, for $u \in B_i$, since $d(u)\geq \frac{n-6}{6}$, we have that $u$ is adjacent to $\Omega(n)$ vertices in $V_{3-i}$.

\begin{claim}\label{CL:Ai} Every vertex in $V_i$ is adjacent to at most one vertices in $A_i$. 
Moreover, if $uv$ is an edge in $A_i$, then no vertex in $B_i$ is adjacent
to $u$ or $v$.
\end{claim}
\begin{proof}
Assume otherwise, without loss of generality, let $u u_1, u u_2$ be 
edges with $u \in V_1$ and $u_1, u_2 \in A_1$. 
Then $|N_2(u)|=\Omega(n)-o(n)=\Omega(n)$.
Since each vertex in $A_1$ has at least $\left|N_2(u)\right|-o(n)$ neighbors in $N_2(u)$,  vertices $u, u_1, u_2$ have $\left|N_2(u)\right|-o(n)\ge 3$ common 
neighbors in $N_2(u)$.
This yields a copy of $\widehat{P_5}$ with center $u$ in $G$, a contradiction. 
The latter of this claim holds  for similar reasons.

\end{proof}

By Claim~\ref{CL:Ai}, every vertex in  $B_i$ is adjacent to at most one vertex of $A_i$. 
Thus, the total number of vertices in $A_i$ that are adjacent to some vertex in $B_i$ is $o(n)$.
By moving these vertices from $A_i$ to $B_i$, we can ensure that there are no edges between $A_i$
and $B_i$, and that $\Delta(G[A_i])\leq 1$, for $i=1,2$. 

\begin{claim}\label{8271}  $e(G[A_1])+e(G[A_2])=\frac{n}{4}-o(n)$.
\end{claim}
\begin{proof}
  Since $\Delta(G[A_i])\leq 1$,  $G[A_i]$ is triangle-free, for $i=1,2$. On the other hand,
  $\vert B\vert=o(n),$  $t(G[B])=o(n^2).$
  Thus we have the number of tirangles with one edge in $G[V_i]$ and a vertex in $V_{3-i}$
  is $\frac{n^2}{8}-o(n^2)$. We have $(e(G[A_1])+e(G[A_2]))\cdot \vert A_2\vert\geq\frac{n^2}{8}-o(n^2)$, 
  and $e(G[A_1])+e(G[A_2])=\frac{n}{4}-o(n)$.
\end{proof}

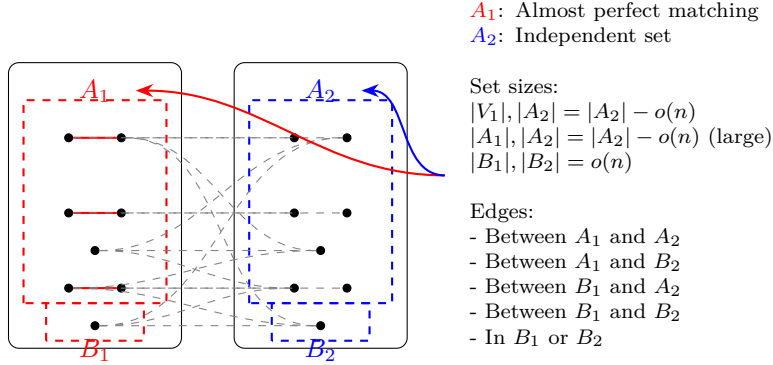
\begin{figure}[http]
\begin{tikzpicture}[
    node distance=1cm,
    vertex/.style={circle, draw, fill=black, inner sep=0pt, minimum size=3pt},
    setlabel/.style={text width=2.5cm, align=left, font=\footnotesize},
    A1/.style={red, thick},
    A2/.style={blue, thick},
    arrow/.style={->, >=Stealth, thick},
    edge/.style={gray, thin, dashed}
]

\foreach \i in {1,2,3} {
    \node[vertex] (a1-left\i) at (0, \i) {};
    \node[vertex] (a1-right\i) at (0.7, \i) {};
}
\foreach \i in {1,2} {
    \node[vertex] (b1\i) at (0.35, \i-0.5) {};
}

\foreach \i in {1,2,3} {
    \node[vertex] (a2-left\i) at (3, \i) {};
    \node[vertex] (a2-right\i) at (3.7, \i) {};
}
\foreach \i in {1,2} {
    \node[vertex] (b2\i) at (3.35, \i-0.5) {};
}

\draw[rounded corners] (-0.8,0.2) rectangle (1.5,4);
\draw[rounded corners] (2.2,0.2) rectangle (4.5,4) coordinate (v2label);

\draw[dashed, A1] (-0.6,0.8) rectangle (1.3,3.5) node[midway,above=1.2cm, red] (a1label) {$A_1$};
\draw[dashed, A1] (-0.3,0.3) rectangle (1.0,0.8) node[midway,below=0.1cm, red] {$B_1$};
\draw[dashed, A2] (2.4,0.8) rectangle (4.3,3.5) node[midway,above=1.2cm, blue] (a2label) {$A_2$};
\draw[dashed, A2] (2.7,0.3) rectangle (4.0,0.8) node[midway,below=0.1cm, blue] {$B_2$};

\foreach \i in {1,2,3} {
    \draw[A1, thick] (a1-left\i) -- (a1-right\i);
}

\foreach \i in {1,2,3} {
    \draw[edge] (a1-left\i) -- (a2-left\i);
    \draw[edge] (a1-right\i) -- (a2-right\i);
}
\foreach \i in {1,2} {
    \draw[edge] (b1\i) -- (b2\i);
}
\foreach \i in {1,2} {
    \draw[edge] (b1\i) to[out=0,in=180] (a2-left1);
    \draw[edge] (b1\i) to[out=0,in=180] (a2-right3);
    
    \draw[edge] (b2\i) to[out=180,in=0] (a1-left1);
    \draw[edge] (b2\i) to[out=180,in=0] (a1-right3);
}

\node[setlabel, right=0.5cm of v2label, anchor=west, yshift=-1.5cm] (descr) {
    \begin{tabular}{l}
        \textcolor{red}{$A_1$}: Almost perfect matching \\
        \textcolor{blue}{$A_2$}: Independent set \\
        \\
        Set sizes: \\
        $|V_1|, |A_2| = \vert A_2\vert-o(n)$ \\
        $|A_1|, |A_2| = \vert A_2\vert-o(n)$ (large) \\
        $|B_1|, |B_2| = o(n)$ \\
        \\
        Edges: \\
        - Between $A_1$ and $A_2$ \\
        - Between $A_1$ and $B_2$ \\
        - Between $B_1$ and $A_2$ \\
        - Between $B_1$ and $B_2$\\
        - In $B_1$ or $B_2$
    \end{tabular}
};

\draw[arrow, red] (descr.west) to[out=180,in=0] ([xshift=0.2cm]a1label.east);
\draw[arrow, blue] (descr.west) to[out=180,in=0] ([xshift=0.2cm]a2label.east);

\end{tikzpicture}
\caption{Structure of $G$}
\label{Fig1}
\end{figure}

By Claim \ref{8271}, without loss of generality, we may assume $\Delta(G[A_1])= 1$, 
$e(G[A_2])=0$, and $e(G[A_1])=\frac{n}{4}-o(n)$. 
Indeed, if {$e(G[V_1])\geq 2$}, let $ab, cd$ be two 
independent edges in $A_1$. Since every vertex in $A_1$ is adjacent to all but 
at most $o(n)$ vertices of $A_2$,  let $D_2$ denote the common
neighborhood of $a, b, c, d$ in $A_2$, then $\left|D_2\right|=\vert A_2\vert-o(n)$. 
{If there 
is an edge in $D_2$, say $ef$, then we could find a copy of $\widehat{P_5}$ with center $e$, a contradiction.} 
Thus there is no edge in $D_2$, and the number of edges in $A_2$ is $o(n)$.
Similarly, if $e(G[A_2])\geq 2$, then the number of edges in $A_1$ is $o(n)$, and we have $e(G[A_1])+e(G[A_2])=o(n).$
{This is contradictory to Claim \ref{8271}.}
Thus there exists $j$,
such that $e(G[A_j])\leq 1$. Assume $j=2$. If $e(G[A_2])= 1$, let $gh$ be an edge
in {$G[A_2]$}, let $D_1$  be the common neighborhood of $g,h$ in $A_1$, then $\vert D_1\vert=\vert A_2\vert-o(n)$.
If there are two edges $kl, mn$ in $D_1$, since $\Delta(G[A_1])\leq 1$, they are independent.
We can find a copy of $\widehat{P_5}$ with center $g$ in $G$, a contradiction.
 Thus there is at most one edge in $D_1$, and the number of edges in $A_1$ is $o(n)$,  $e(G[A_1])+e(G[A_2])=o(n)$
 which is contradictory to Claim \ref{8271}. Thus we may assume $A_2$ is a stable 
 set in $G$, and $e(G[A_1])=\frac{n}{4}-o(n)$. The structure of graph G is shown in Figure \ref{Fig1}.

\begin{claim}\label{CL:Bi} 
For  each vertex $u\in B_1$, $t(u)\leq  d_{2}(u)+o(n).$
\end{claim}
\begin{proof} 
Frist, we would show $t_1(u)\leq  d_{2}(u)+o(n).$
We proceed by induction on $d_1(u)$.
If $d_1(u)=1$, let $N_1(u)=\{x\}$,
then $t_1(u)\leq d_2(u).$


Now assume $d_1(u)\geq 2$ and the claim is true for vertex with degree in $B_1$  smaller than $d_1(u)$. 
{Since $G$ is $\widehat{P_5}$-free, any two vertices $x$ and $y$ in $N_1(u)$ have at most two common neighbors in $N_2(u)$.}

\noindent \textbf{Case 1:} There exist two vertices $x$ and $y$ in $N_1(u)$ such that 
$\vert N_2(u)\cap N_2(x)\cap N_2(y)\vert=2$. 

Then neither $x$ nor $y$ can have other
neighbors in $N_2(u)$. Let $\{x', y'\}=N_2(u)\cap N_2(x)\cap N_2(y)$. 
For any vertex $z\in N_1(u)\setminus \{x,y\}$, we have $zx', zy'\notin E(G)$, since
$G$ is $\widehat{P_5}$-free. Let $G'=G-\{x,y,x', y'\}$.  By 
the induction hypothesis on $u$ in $G'$, we have
$$t_1(u)-t(xu)-t(yu)\leq  d_{2}(u)-2+o(n).$$
Combining with $t(xu)+t(yu)\leq 4+o(n)$,
we have $$t_1(u)\leq  d_{2}(u)+o(n).$$

\noindent \textbf{Case 2:} There exist two vertices $x$ and $y$ in $N_1(u)$ such that 
$\vert N_2(u)\cap N_2(x)\cap N_2(y)\vert=1$. 

Denote $N_2(u)\cap N_2(x)\cap N_2(y)=\{x'\}$.
Let $X=N_2(x)\cap N_2(u)\setminus\{x^\prime\}$ and $Y=N_2(y)\cap N_2(u)\setminus\{x^\prime\}$. For any vertex $z\in N_1(u)\setminus \{x,y\}$, we must have 
$N_2(z)\cap X=\emptyset$ and $N_2(z)\cap Y=\emptyset$; otherwise we can find a copy of 
$\widehat{P_5}$ in $G$, a contradiction. Similarly, remove $x,y,X, Y$ from the graph and
 use the induction hypothesis, we have 
$$t_1(u)-t(xu)-t(yu)\leq { d_{2}(u)-\vert X\vert-\vert Y\vert +o(n).}$$
Combining with $t(xu)+t(yu)\leq \vert X\vert +\vert Y\vert +2 +o(n)$,
we have $$t_1(u)\leq  d_{2}(u)+o(n).$$

\noindent \textbf{Case 3:} For any vertices $x,y\in N_1(u)$, they are no common neighbors in $N_2(u)$.
Obviously,  
$t_1(u)\leq d_2(u)$.

Next, we would show $t(u)\leq  d_{2}(u)+o(n)$.
Since $A_2$ is a stable set in $G$, the second or third type triangles of $u$ must contain an edge in $G[N_2(u)\cap B_2]$ or in $G[N_2(u)\cap B_1]$.
Since $\vert B_i\vert =o(n)$, 
$t_2(u)+t_3(u)\leq \frac{3}{2}\vert B_2\vert+\frac{3}{2}\vert B_1\vert=o(n)$, so $t(u)\leq  d_{2}(u)+o(n)$.
\end{proof}

\begin{claim}\label{CL:Ht} 
Let $H$ be an induced  subgraph of $G[B_1]$.
Let $t^\prime(H)$ be the number of triangles $T$ in $G$ of the following two types:

1. $\vert V(T)\cap V(H)\vert \geq 2$,

2. $\vert V(T)\cap V(H)=1\vert$, $\vert V(T)\cap B_2\vert =2$.

If $\vert V(H) \vert=2k+1$,
then $t(H)\leq k\cdot \vert A_2\vert+o(n)$; 
if $\vert V(H) \vert=2k+2$, then  $t(H)\leq (k+1)\cdot \vert A_2\vert-\Omega(n)$.
\end{claim}

\noindent\textbf{Remark.} If $H$ is a connected component in $G[B_1]$, then $t^\prime(H)=t(H)$.
If not, then $t^\prime(H)<t(H)$ may hold.

\begin{proof}
  First we would show that if the statement is true for the odd case, 
 it is also true for the even case. 
  By Claim~\ref{CL:Ai}, vertices in $B_1$ with $r(u)=o(n)$( these vertices are in $A_1$ initially)
  can not be adjacent. We may assume there is a vertex $u$ with 
  $r(u)=\Omega(n)$ in  the selected $2k+2$ vertices since otherwise $H$ is an empty graph, and
  $t^\prime(H)\leq (2k+2)\cdot \frac{3}{2}\vert B_2\vert \leq (k+1)\cdot \vert A_2\vert-\Omega(n)$.

Let $H^\prime=G[V(H)\setminus\{u\}]$. Then $|V(H')|=2k+1$ and  $t^\prime(H')\leq k\cdot \vert A_2\vert+o(n)$. By Claim~\ref{CL:Bi}, 
$$t(u)\leq  d_{2}(u)+o(n){\leq |A_2|-r(u)+o(n)}.$$
Combining these two inequalities,
we have  $t(H)\leq (k+1)\cdot \vert A_2\vert-\Omega(n)$.

Next, we establish the statement for the odd case by induction on $k$. 
Let $K$ be the set of  the selected $2k+1$ vertices. 
The base case when $k=0$ is trivial. Since we may assume $K=\{u\}$, and
$t^\prime(H)=e(G[N_2(u)\cap B_2])=o(n)$. For $k\ge 1$, suppose the claim is true for $2k^\prime+1$ with  $k^\prime<k$. 
If there exist  two adjacent vertices $x,y$ in $K$ such that 
$d_2(x)+d_2(y)\leq \vert A_2\vert+o(n)$, then by Claim~\ref{CL:Bi}, 
$$t(x)+t(y)\leq  d_{2}(x)+d_2(y)+o(n)\leq \vert A_2\vert+o(n).$$
Using the induction hypothesis on $K\setminus\{x,y\}$,
we have $$t^\prime(G[K\setminus\{x,y\}])\leq (k-1)\cdot \vert A_2\vert+o(n).$$
Combining these two  inequalities, we have the desired claim.
Thus we may assume for any two adjacent vertices $x,y$ in $H$, $d_2(x)+d_2(y)\geq\vert A_2\vert+\Omega(n)$. This yields that $x$ and $y$ have $\Omega(n)$ common neighbors in $A_2$.

In addition, if there is a vertex $u\in K$, such that $\vert N_H(u)\vert =1$,
remove $u$ and its neighbor $v$ from $H$ and use the hypothesis induction, we have 
 $$t^\prime(G[K\setminus\{u,v\}])\leq (k-1)\cdot \vert A_2\vert+o(n).$$
By Claim~\ref{CL:Bi}, $t(v)\leq  d_{2}(v)+o(n)\leq \vert A_2\vert+o(n).$
Since $t_2(u)=o(n)$, we have
 $t^\prime(H)\leq t^\prime(G[K\setminus\{u,v\}])+t(v)+t_2(u) \leq k\cdot \vert A_2\vert+o(n)$.
Therefore, we may further assume that for any vertex $u\in K$, $\vert N_H(u)\vert \geq 2$. 
Choose $v,w\in N_H(u)$. Since $u$ and $v$ (resp. $w$) have $\Omega(n)$ common neighbors in $A_2$, we have
$$N_2(u)\cap N_2(v)\cap N_2(w)= \emptyset;$$
otherwise, we can find a copy of $\widehat{P_5}$ with center $u$ in $G$, a contradiction.
Namely, any two vertices in $N_1(u)$ cannot have common neighbors in $N_2(u)$. 
In other words, any vertex in $N_2(u)$ is adjacent to at most one vertex of $N_1(u)$ (as shown in Fig.~\ref{Fig:N2u}). 
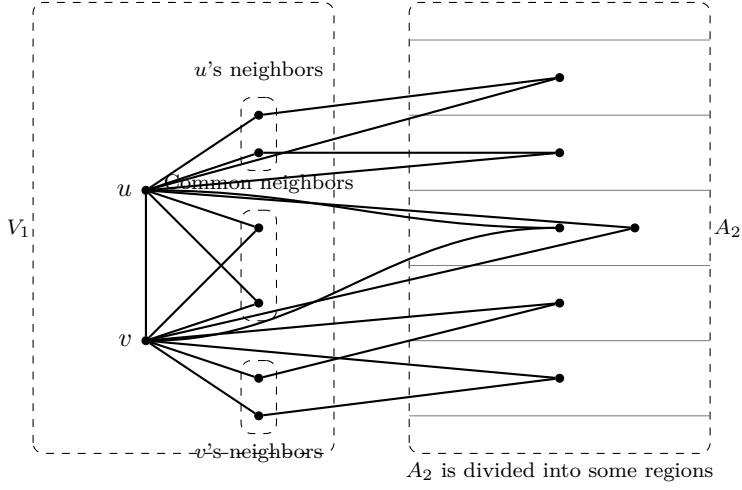
\begin{figure}[http]
\begin{tikzpicture}[
    node distance=0.5cm,
    vertex/.style={circle, draw, fill=black, inner sep=0pt, minimum size=3pt},
    setlabel/.style={text width=2.5cm, align=left, font=\footnotesize},
    edge/.style={black, thick},
    neighborbox/.style={draw, dashed, rounded corners, inner sep=5pt},
    region/.style={draw, gray, thin}
]

\node[setlabel] at (-1.6, 2.5) {$V_1$};
\node[setlabel] at (7.8, 2.5) {$A_2$};

\draw[rounded corners, dashed] (-2.5,-0.5) rectangle (1.5,5.5); 
\draw[rounded corners, dashed] (2.5,-0.5) rectangle (6.5,5.5);  

\node[vertex, label=left:$u$] (u) at (-1,3) {};
\node[vertex, label=left:$v$] (v) at (-1,1) {};
\draw[edge] (u) -- (v); 

\node[vertex] (u1) at (0.5,4) {};
\node[vertex] (u2) at (0.5,3.5) {};
\draw[edge] (u) -- (u1);
\draw[edge] (u) -- (u2);

\node[vertex] (v1) at (0.5,0.5) {};
\node[vertex] (v2) at (0.5,0) {};
\draw[edge] (v) -- (v1);
\draw[edge] (v) -- (v2);

\node[vertex] (common1) at (0.5,2.5) {};
\node[vertex] (common2) at (0.5,1.5) {};
\draw[edge] (u) -- (common1);
\draw[edge] (u) -- (common2);
\draw[edge] (v) -- (common1);
\draw[edge] (v) -- (common2);

\node[neighborbox, fit=(u1) (u2)] (ubox) {};
\node[above=0.1cm of ubox] {\footnotesize $u$'s neighbors};

\node[neighborbox, fit=(v1) (v2)] (vbox) {};
\node[below=0.0001cm of vbox] {\footnotesize $v$'s neighbors};

\node[neighborbox, fit=(common1) (common2)] (cbox) {};
\node[above=0.1cm of cbox] {\footnotesize Common neighbors};

\foreach \i in {0,1,2,3,4,5} {
    \draw[region] (2.5,\i) -- (6.5,\i); 
}

\node[vertex] (w1) at (4.5,0.5) {}; 
\node[vertex] (w2) at (4.5,1.5) {}; 
\node[vertex] (w3) at (4.5,2.5) {}; 
\node[vertex] (w4) at (4.5,3.5) {}; 
\node[vertex] (w5) at (4.5,4.5) {}; 
\node[vertex] (w6) at (5.5,2.5) {}; 

\draw[edge] (u) -- (w4);
\draw[edge] (u2) -- (w4);
\draw[edge] (u1) -- (w5);
\draw[edge] (u) -- (w5);
\draw[edge] (u) -- (w6); 

\draw[edge] (v) -- (w1);
\draw[edge] (v1) -- (w2);
\draw[edge] (v2) -- (w1);
\draw[edge] (v) -- (w2);
\draw[edge] (v) -- (w6); 

\draw[edge] (u) to[out=0,in=180] (w3);
\draw[edge] (v) to[out=0,in=180] (w3);


\node[below] at (4.5,-0.5) {\footnotesize $A_2$ is divided into some regions};

\end{tikzpicture}
\caption{Structure between $N_1(u)$ and $N_2(u)$}
\label{Fig:N2u}
\end{figure}

Therefore,
$$\sum_{x\in N_1(u),x\neq v}t(xu)\leq \vert N_2(u)\setminus N_2(v)\vert+o(n),$$
$$\sum_{x\in N_1(v),x\neq u}t(xv)\leq \vert N_2(v)\setminus N_2(u)\vert+o(n).$$
Consequently, 
\[
\begin{aligned}
t(G[\{u,v\}])&\leq \sum_{\substack{x \in N_1(u) \\ x \neq v}} t(xu) 
+ \sum_{\substack{x \in N_1(v) \\ x \neq u}} t(xv) 
+ t(uv)+o(n)\\
&\leq \lvert N_2(u) \setminus N_2(v) \rvert 
+ \lvert N_2(v) \setminus N_2(u) \rvert 
+ \lvert N_2(u) \cap N_2(v) \rvert+o(n) \\
&\leq \vert A_2\vert+o(n)
\end{aligned}
\]
Removing the vertices $\{u,v\}$ from $K$ and applying the induction hypothesis to  $K\setminus\{u,v\}$,
we have $t^\prime(G[K\setminus\{u,v\}]) \leq (k-1)\cdot \vert A_2\vert+o(n)$.
Combining the inequalities above, we have
$t^\prime(H)\leq k\cdot\vert A_2\vert+o(n)$.
\end{proof}

\begin{claim}\label{CL:B1=empty} 
	$B_1=\emptyset$.
\end{claim}
\begin{proof}
  Suppose not. Let $H$ be a connected  component in $G[B_1]$. 
  If the order of $H$ is odd, say $2k+1$, then we do the following
  adjustment to $G$ to get a graph with more triangles while maintaining $\widehat{P_5}$-free.
  Remove all the edges adjacent to $H$.
  Move a vertex $u$ from $H$ to $A_2$ and make it adjacent to all vertices in $A_1$,
  while let the remaining $2k$ vertices in $H$ form a perfect matching consisting of 
  $k$ edges by pairing them into $k$ disjoint pairs. And make each vertex of the $2k$ vertices  adjacent to all vertices in $A_2$. 
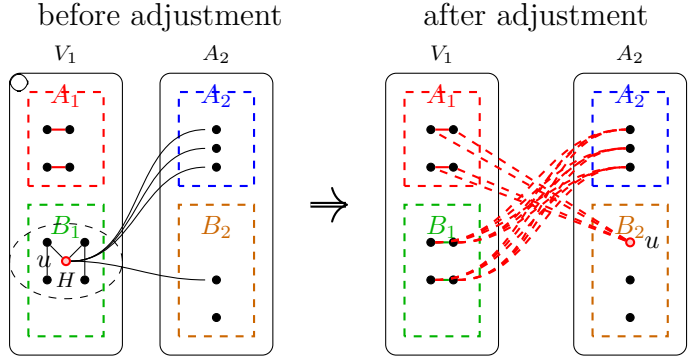
\begin{figure}
\begin{tikzpicture}[
    scale=0.5,
    node distance=0.8cm,
    vertex/.style={circle, draw, fill=black, inner sep=0pt, minimum size=3pt},
    setlabel/.style={text width=2cm, align=center, font=\footnotesize},
    A1/.style={red, thick},
    A2/.style={blue, thick},
    B1/.style={green!70!black, thick},
    B2/.style={orange!80!black, thick},
    edge/.style={black, thin},
    newedge/.style={red, thick, dashed},
    highlight/.style={fill=red!20, draw=red, thick},
    cloud/.style={draw, ellipse, minimum width=1.5cm, minimum height=1cm, dashed}
]

\node[font=\large] at (-1, 9) { before adjustment };

\node[setlabel] at (-3.5, 8) {$V_1$};
\node[setlabel] at (0.5, 8) {$A_2$};

\draw[rounded corners] (-5,0) rectangle (-2,7.5);
\draw[rounded corners] (-1,0) rectangle (2,7.5);

\draw[dashed, A1] (-4.5,4.5) rectangle (-2.5,7) node[midway,above=0.3cm, red] {$A_1$};
\draw[dashed, B1] (-4.5,0.5) rectangle (-2.5,4) node[midway,above=0.3cm, green!70!black] {$B_1$};
\draw[dashed, A2] (-0.5,4.5) rectangle (1.5,7) node[midway,above=0.3cm, blue] {$A_2$};
\draw[dashed, B2] (-0.5,0.5) rectangle (1.5,4) node[midway,above=0.3cm, orange!80!black] {$B_2$};

\node[vertex] (a1-1) at (-3.4, 6) {};
\node[vertex] (a1-2) at (-4, 6) {};
\node[vertex] (a1-3) at (-3.4, 5) {};
\node[vertex] (a1-4) at (-4, 5) {};
\draw[A1] (a1-1) -- (a1-2);
\draw[A1] (a1-3) -- (a1-4);

\node[vertex, highlight, label=left:{$u$}] (u-before) at (-3.5, 2.5) {};
\node[cloud, fit=(u-before)] (H-before) {};
\node[setlabel] at (-3.5, 2) {$H$};
\node[vertex] (h1) at (-3, 3) {};
\node[vertex] (h2) at (-4, 3) {};
\node[vertex] (h3) at (-3, 2) {};
\node[vertex] (h4) at (-4, 2) {};
\draw[edge] (u-before) -- (h1);
\draw[edge] (u-before) -- (h2);
\draw[edge] (h1) -- (h3);
\draw[edge] (h2) -- (h4);

\foreach \i in {1,2,3} {
    \node[vertex] (a2-\i) at (0.5, 6.5-\i*0.5) {};
}

\node[vertex] (b2-after-1) at (0.5, 2) {};
\node[vertex] (b2-after-2) at (0.5, 1) {};

\draw[edge] (u-before) to[out=0, in=180] (0.2, 6);
\draw[edge] (u-before) to[out=0, in=180] (0.2, 5.5);
\draw[edge] (u-before) to[out=0, in=180] (0.2, 5);
\draw[edge] (u-before) to[out=0, in=180] (0.2, 2);

\draw[->, thick,double, black] (3,4) -- (4,4);
\node[font=\large] at (9, 9) {after adjustment};

\node[setlabel] at (6.5, 8) {$V_1$};
\node[setlabel] at (11.5, 8) {$A_2$};

\draw[rounded corners] (5,0) rectangle (8,7.5);
\draw[rounded corners] (10,0) rectangle (13,7.5);

\draw[dashed, A1] (5.5,4.5) rectangle (7.5,7) node[midway,above=0.3cm, red] {$A_1$};
\draw[dashed, B1] (5.5,0.5) rectangle (7.5,4) node[midway,above=0.3cm, green!70!black] {$B_1$};
\draw[dashed, A2] (10.5,4.5) rectangle (12.5,7) node[midway,above=0.3cm, blue] {$A_2$};
\draw[dashed, B2] (10.5,0.5) rectangle (12.5,4) node[midway,above=0.3cm, orange!80!black] {$B_2$};

\node[vertex] (a1-after-1) at (6.2, 6) {};
\node[vertex] (a1-after-2) at (6.8, 6) {};
\node[vertex] (a1-after-3) at (6.2, 5) {};
\node[vertex] (a1-after-4) at (6.8, 5) {};
\draw[A1] (a1-after-1) -- (a1-after-2);
\draw[A1] (a1-after-3) -- (a1-after-4);

\node[vertex] (b1-1) at (6.2, 3) {};
\node[vertex] (b1-2) at (6.8, 3) {};
\node[vertex] (b1-3) at (6.2, 2) {};
\node[vertex] (b1-4) at (6.8, 2) {};
\draw[B1, thick] (b1-1) -- (b1-2);
\draw[B1, thick] (b1-3) -- (b1-4);

\foreach \i in {1,2,3} {
    \node[vertex] (a2-after-\i) at (11.5, 6.5-\i*0.5) {};
}

\node[vertex, highlight, label=right:{$u$}] (u-after) at (11.5, 3) {};
\node[vertex] (b2-after-1) at (11.5, 2) {};
\node[vertex] (b2-after-2) at (11.5, 1) {};

\foreach \i in {1,...,4} {
    \draw[newedge] (u-after) -- (a1-after-\i);
}

\foreach \i in {1,3} {
    \foreach \j in {1,2,3} {
        \draw[newedge] (b1-\i) to[out=0, in=180] (a2-after-\j);
        \draw[newedge] (b1-\the\numexpr\i+1\relax) to[out=0, in=180] (a2-after-\j);
    }
}


\node[draw, rounded corners, fill=white, fill opacity=0.8, text opacity=1, 
      align=left, font=\scriptsize, anchor=north west] at (-5, 7.5) {
};

\end{tikzpicture}
\caption{Structural adjustment diagram}
\label{Fig3}
\end{figure}

  Obviously, after this adjustment,   no new copy of $\widehat{P_5}$  will be created.
 By Claim~\ref{CL:Ht}, before adjustment,  $t(H)\leq k\cdot \vert A_2\vert+o(n)$.
  After adjustment, 
  $t(H)=\frac{n}{4}-o(n)+k\cdot\vert A_2\vert$, since there are $\frac{n}{4}-o(n)$ edges in $A_1$.
  This contradicts that $t(G)=\ex(n, K_3, \widehat{P_5})$.

  Similarly,  if the order of $H$ is even, say $2k$, then we do the following
  adjustment to $G$ to get a graph with more triangles while maintaining $\widehat{P_5}$-free.
  Remove all edges adjacent to $H$, add a perfect matching to $H$, and make every vertex in $H$ adjacent to every vertex in $A_2$.
Then by Claim~\ref{CL:Ht}, before adjustment $t(H)\leq k\cdot \vert A_2\vert-\Omega(n)$,
and after adjustment $t(H)= k\cdot\vert A_2\vert$. Thus we
get a graph with more triangles while maintaining $\widehat{P_5}$-free, which contradicts that $t(G)=\ex(n, K_3, \widehat{P_5})$.
\end{proof}

By Claim~\ref{CL:B1=empty}, for any vertex $u\in V_1$, $d_1(u)=1$, namely there 
is a perfect matching in $G[V_1]$. Since otherwise, if there is an isolated vertex in 
 $G[V_1]$, we can move it to $A_2$, and make it adjacent to each vertex of $V_1$, then 
 we get a graph with more triangles while maintaining $\widehat{P_5}$-free.

 In addition, for each vertex $w\in B_2$, $r(w)=\Omega(n)$. Indeed, suppose $xy$ is 
 an edge in $G[B_2]$ and $r(x)=o(n)$, then $\vert N_1(x)\vert =\vert A_2\vert-o(n)$, 
 and $e(G[N_1(x)])=\frac{n}{4}-o(n)$. Since $d_1(y)=\Omega(n)$, we can easily find 
 two independent edges $z_1z_2, z_3z_4$ in $G[N_1(x)]$ such that $z_1, z_3\sim y$. 
 Then $z_2z_1yz_3z_4$ is a copy of $P_5$ in $N(x)$, a contradiction.

\begin{claim}
$B_2=\emptyset$.
\end{claim}
\begin{proof}
Otherwise, we do the following
adjustment to $G$ to get a graph with more triangles while maintaining $\widehat{P_5}$-free.
Remove all the edges in $B_2$, add edges to make each vertex of $B_2$ adjacent to each vertex of $V_1$.
Obviously, after the adjustment, no copy of $\widehat{P_5}$  will be created. By Lemma \ref{8291}, before adjustment,
$t(B_2)< \vert B_2\vert \cdot \frac{\vert V_1\vert}{2}.$
After adjustment, 
$t(B_2)= \vert B_2\vert \cdot \frac{\vert V_1\vert}{2}.$
 This contradicts that $t(G)=ex(n, K_3, \widehat{P_5})$. Thus $B_2=\emptyset$.
\end{proof}

At this point, we have completed the structural characterization of graph $G$:
$V(G)=V_1\cup A_2$, $\vert V_1\vert=\vert A_2\vert-o(n)$, $\vert A_2\vert=\vert A_2\vert-o(n)$, $G[V_1]$ has a perfect matching, and $G[A_2]$
is a stable set in $G$. It follows from a straightforward calculation that 
$t(G)\leq \left\lfloor \frac{n^2}{8} \right\rfloor$, and the equality holds if and only if $G\cong H_n$.

\section{Remarks and Discussions}

In this paper, using the stability method, we show that for sufficiently large $n$, $\ex(n,K_3, \widehat{P_5})=\left\lfloor n^2/8\right\rfloor .$ 
and the extremal graph is unique. 
Determining the exact value of $\ex(n, K_3, \widehat{P_k})$ for $k\geq 6$
 is a noteworthy and challenging problem.
We also believe that the  lower bound in Proposition \ref{lemma 1} is asymptotically tight for all fixed $ k\geq 6$.
And we find a better lower bound than Construction~\ref{CONS:Fnk}.
\begin{constr}\label{8292} Let $H_{n, k}$ be the graph formed by  the complete bipartite graph $K_{n/2, n/2}$ with partition $(A, B)$,  
together with additional edges in $A$ such that $H_{n, k}[A]$ consists of disjoint copies of $K_{\left\lfloor\frac{k}{2}\right\rfloor}$,  
where $n$ is a multiple of $2\left\lfloor \frac{k}{2}\right\rfloor$.
\end{constr}	
Obviously,  $H_{n, k}$ is $\widehat{P}_k$-free, since the neighborhood of every vertex in $B$ is a disjoint union of $K_{\left\lfloor\frac{k}{2}\right\rfloor}$, 
and the neighborhood of every vertex in $A$ is isomorphic to $K_{\left\lfloor\frac{k}{2}\right\rfloor-1, \vert A_2\vert}$.

From Construction~\ref{8292}, we have 
$$
e(A)=\frac{n}{4\left\lfloor \frac{k}{2}\right\rfloor}\cdot \left\lfloor \frac{k}{2}\right\rfloor\left(\left\lfloor \frac{k}{2}\right\rfloor-1\right).$$
Since every triangle of $H_{n, k}$ consists of an edge from $H_{n, k}[A]$ and a vertex from $B$, or is entirly in $A$,  we have 
$$
t\left(H_{n, k}\right)=\vert A_2\vert\cdot \frac{n}{4\left\lfloor \frac{k}{2}\right\rfloor}\cdot \left\lfloor \frac{k}{2}\right\rfloor\left(\left\lfloor \frac{k}{2}\right\rfloor-1\right)+
\frac{n}{2\left\lfloor \frac{k}{2}\right\rfloor}\cdot\binom{\left\lfloor \frac{k}{2}\right\rfloor}{3}>t(\left(F_{n, k}\right)).
$$

\vspace{5pt}
\noindent{\bf Acknowledgements}:
This work was supported by the National Key Research and Development Program of China (2023YFA1010203), the National Natural Science Foundation of China (No. 12401455, 12471336), and the Innovation Program for Quantum Science and Technology (2021ZD0302902). 


\bibliographystyle{plain}  
\bibliography{ref}  

\end{document}